\documentclass{amsart}
\usepackage{latexsym}

\usepackage{amsfonts,amssymb,amsmath}
\usepackage{amsthm}
\usepackage{graphicx,color}
\usepackage{color}
\usepackage{amsmath} 
\newtheorem{theorem}{Theorem}	
\newtheorem{lemma}{Lemma}[section]
\newtheorem{corollary}{Corollary}		
\newtheorem{proposition}{Proposition}

\newtheorem{definition}{Definition}

\newfont{\bg}{cmr9 scaled\magstep2}
\newcommand{\bigzerol}{\smash{\lower1.0ex\hbox{\bg 0}}}

\newcommand{\R}{\mathbb{R}}

\makeatletter
\newcommand{\tpitchfork}{%
  \vbox{
    \baselineskip\z@skip
    \lineskip-.52ex
    \lineskiplimit\maxdimen
    \m@th
    \ialign{##\crcr\hidewidth\smash{$-$}\hidewidth\crcr$\pitchfork$\crcr}
  }%
}
\makeatother
\title[Preservation of {immersed or} injective properties]
{
Preservation of {immersed} or injective properties by composing generic 
generalized distance-squared mappings
}
\author{Shunsuke Ichiki}
\thanks{The first author is Research Fellow DC1 of Japan Society for the Promotion of Science }
\address{
Graduate School of Environment and Information Sciences,  
Yokohama National University, 
Yokohama 240-8501, Japan}
\email{ichiki-shunsuke-jb@ynu.jp}

\author{Takashi Nishimura}
\address{Research Group of Mathematical Sciences,
Research Institute of Environment and Information Sciences,
Yokohama National University,
Yokohama 240-8501, Japan}
\email{nishimura-takashi-yx@ynu.jp}

\begin{document}
\date{}
\begin{abstract}
Any generalized distance-squared mapping of equidimensional 
case has singularities, and {\color{black}their} singularity 
type{\color{black}s are} wrapped into mystery
in higher dimensional cases. 
Any generalized distance-squared mapping of equidimensional 
case is not injective. 
Nevertheless, in this paper, it is shown
 that the {\color{black}non-singular} property or the injective property of a
 mapping is preserved 
 by composing a generic generalized distance-squared mapping of
 equidimensional case. 
\end{abstract}
\maketitle
\noindent
\section{Introduction}\label{section 1}
Throughout this paper, $i$, $j$, $\ell$, $m$ 
{\color{black}and} $n$ stand for positive integers. 
In this paper, unless otherwise stated, all manifolds and mappings belong to class $C^{\infty}$ 
and all manifolds are without boundary. 
Let $p_i=(p_{i1}, p_{i2}, \ldots, p_{im})$  $(1\le i\le \ell)$ 
(resp., $A=(a_{ij})_{1\le i\le \ell, 1\le j\le m}$) 
be a point of $\mathbb{R}^m$ 
(resp., an $\ell\times m$ matrix with non-zero entries).
Set $p=(p_1,p_2,\ldots,p_{\ell})\in (\mathbb{R}^m)^{\ell}$. 
Let $G_{(p, A)}:\mathbb{R}^m \to \mathbb{R}^\ell$ be the mapping 
defined by 
{\small 
\[
G_{(p, A)}(x)=\left(
\sum_{j=1}^m a_{1j}(x_j-p_{1j})^2, 
\sum_{j=1}^m a_{2j}(x_j-p_{2j})^2, 
\ldots, 
\sum_{j=1}^m a_{\ell j}(x_j-p_{\ell j})^2
\right), 
\] }where $x=(x_1, x_2, \ldots, x_m)\in \mathbb{R}^m$. 
The mapping $G_{(p, A)}$ is called a {\it generalized distance-squared mapping}, 
and the $\ell$-tuple of points $p=(p_1,\ldots ,p_{\ell})\in (\mathbb{R}^m)^{\ell}$ 
is called the {\it central point} 
of the generalized distance-squared mapping $G_{(p,A)}$. 
A {\it distance-squared mapping} $D_p$ 
(resp., {\it Lorentzian distance-squared mapping}
$L_p$) is the mapping $G_{(p,A)}$ 
satisfying that each entry of $A$ is $1$ 
(resp., $a_{i1}=-1$ and $a_{ij}=1$ $(j\ne 1)$). 

In \cite{D} (resp., \cite{L}), 
a classification result on distance-squared mappings $D_p$ 
(resp., Lorentzian distance-squared mappings
$L_p$) is given. 

In \cite{G1}, 
a classification result on 
generalized distance-squared mappings of the plane into the plane 
is given. 
If the rank of $A$ is two,  
a generalized distance-squared mapping 
having a generic central point 
is a mapping of which any singular point is 
a {\color{black}{\rm fold point}} 
except one {\color{black}{\rm cusp point}}.     
The singular set is a rectangular hyperbola. 
If the rank of $A$ is one, 
a generalized distance-squared mapping 
having a generic central point is $\mathcal{A}$-equivalent 
to 
{\color{black}{\rm the normal form of fold singularity}} $(x_1,x_2)
\mapsto (x_1,x_2^2)$. 

In \cite{G2}, 
a classification result on 
generalized distance-squared mappings 
of $\mathbb{R}^{m+1}$ into $\mathbb{R}^{2m+1}$ 
is given. 
If the rank of $A$ is $m+1$, 
a generalized distance-squared mapping 
having a generic central point 
is $\mathcal{A}$-equivalent to 
{\color{black}{\rm the normal form of Whitney umbrella}} 
$
(x_1,\ldots ,x_{m+1})\mapsto (x_1^2,x_1x_2,\ldots ,x_1x_{m+1},x_2,\ldots ,x_{m+1})
$.
If the rank of $A$ is less than $m+1$, 
a generalized distance-squared mapping 
having a generic central point 
is $\mathcal{A}$-equivalent to the inclusion 
$
(x_1,\ldots ,x_{m+1})\mapsto (x_1,\ldots ,x_{m+1},0,\ldots ,0)
$. 
\par 
In \cite{G2} and \cite{G1}, 
the properties of {\color{black}generic} generalized distance-squared mappings 
are investigated. 
Hence, it is natural to investigate 
the properties of compositions {\color{black}with} generic generalized distance-squared mappings 

We have another original motivation. 
Height functions and distance-squared functions have been investigated 
in detail {\color{black}so far}, 
and they are a useful tool 
in the applications of singularity theory to differential geometry 
({\color{black}for instance}, see \cite{CS}). 
The mapping in which each component is a height function is 
{\color{black}nothing but} a projection. 
In \cite{GP}, compositions of 
{\color{black}generic} projections and embeddings are investigated.

On the other hand, 
the mapping in which each component 
is a distance-squared function is a distance-squared mapping. 
{\color{black}And}, the notion of generalized distance-squared mapping is 
an extension of the distance-squared mappings. 
Therefore, 
{\color{black}it is again natural to} investigate compositions 
{\color{black}with generic} generalized distance-squared mappings. 

Any generalized distance-squared mapping of equidimensional 
case $G_{(p,A)}:\mathbb{R}^m\to \mathbb{R}^m$ has singularities (see Lemma \ref{singularities} in Appendix).
 Nevertheless, in Theorem \ref{immersion}, it is shown
 that the {\color{black}immersed} property of a
 mapping is preserved by composing a generic generalized distance-squared mapping of
 equidimensional case. 

\begin{theorem}\label{immersion}
Let $N$ be an $n$-dimensional manifold, 
and let $f:N\to \R^m$ be an immersion $(m\geq 2n)$. 
Then, there exists a subset $\Sigma$ of $(\R^{m})^m$ with 
Lebesgue measure zero such that for any $p\in (\R^{m})^m -\Sigma $, 
the composition $G_{(p,A)}\circ f:N\to \R^m$ is an immersion. 
\end{theorem}

Any generalized distance-squared mapping of equidimensional 
case $G_{(p,A)}:\mathbb{R}^m\to \mathbb{R}^m$ is not injective (see Lemma \ref{not injective} in Appendix). 
Nevertheless, in Theorem \ref{injective}, it is shown
 that the injective property of a
 mapping is preserved by composing a generic generalized distance-squared mapping of equidimensional case.

\begin{theorem}\label{injective}
Let $N$ be an $n$-dimensional manifold, 
and let $f:N\to \R^m$ be injective $(m\geq 2n+1)$. 
Then, there exists a subset $\Sigma$ of $(\R^m)^m$ with 
Lebesgue measure zero such that for any $p\in (\R^m)^m -\Sigma $, 
the composition $G_{(p,A)}\circ f:N\to \R^m$ is injective. 
\end{theorem}

By combining Theorem \ref{immersion} and Theorem \ref{injective}, 
we have the following proposition. 
\begin{proposition}\label{injective immersion}
Let $N$ be an $n$-dimensional manifold, 
and let $f:N\to \R^m$ be an injective immersion $(m\geq 2n+1)$. 
Then, there exists a subset $\Sigma$ of $(\R^m)^m$ with 
Lebesgue measure zero such that for any $p\in (\R^m)^m -\Sigma $, 
the composition $G_{(p,A)}\circ f:N\to \R^m$ is an injective immersion. 
\end{proposition}

\subsection{Remark}
Suppose that the mapping $G_{(p,A)}\circ f:N\to \mathbb{R}^m$ is proper 
in Proposition \ref{injective immersion}. 
Then, the injective immersion of $G_{(p,A)}\circ f$ implies the embedding of it 
(see \cite{GG}, p.11). 
Hence, we have the following as a corollary of Proposition \ref{injective immersion}. 
\begin{corollary}\label{injective immersion}
Let $N$ be an $n$-dimensional compact manifold, 
and let $f:N\to \R^m$ be an embedding $(m\geq 2n+1)$. 
Then, there exists a subset $\Sigma$ of $(\R^m)^m$ with 
Lebesgue measure zero such that for any $p\in (\R^m)^m -\Sigma $, 
the composition $G_{(p,A)}\circ f:N\to \R^m$ is an embedding. 
\end{corollary}
\par 
\bigskip 
In Section \ref{section 2}, it is reviewed some of standard definitions, 
and an important lemma for the proofs of Theorem \ref{immersion} and 
Theorem \ref{injective} is given. 
Section \ref{section 3} (resp., Section \ref{section 4}) devotes 
the proof of Theorem \ref{immersion} (resp., Theorem \ref{injective}). 
Finally, in Subsection \ref{singularities} (resp., Subsection \ref{not injective}), 
{\color{black}for the sake of readers' convenience, it is given}
the proof that 
any generalized distance-squared mapping of equidimensional case 
has singularities (resp., the proof that 
any generalized distance-squared mapping of equidimensional case 
is not injective).   
\section{Preliminaries}\label{section 2}
Let $N$ and $P$ be manifolds and let $J^r(N,P)$ be the space of 
$r$-jets of mappings of $N$ into $P$. 
For a given mapping $g:N\to P$, 
the mapping $j^r g:N\to J^r(N,P)$ 
is defined by $q \mapsto j^r g(q)$  
(for details on the space $J^r(N,P)$ or 
the mapping $j^r g:N\to J^r(N,P)$, see for example \cite{GG}). 

Next, we recall the definition of transversality. 
\begin{definition}\label{transverse}
{\rm Let $W$ be a submanifold of $P$. 
For a given mapping $g:N\to P$, 
we say that $g:N\to P$ is {\it transverse} to $W$ 
if for any $q\in N$, $g(q)\not\in W$ or  
in the case of $g(q)\in W$, the following holds: 
\begin{eqnarray*}
dg_q(T_qN)+T_{g(q)}W=T_{g(q)}P.
\end{eqnarray*}
}
\end{definition}
For the proofs of Theorem \ref {immersion} and 
Theorem \ref{injective}, the following lemma is important. 
\begin{lemma}[\cite{abra}, \cite{GP}]\label{abra}
Let $N$, $P$, $Z$ be manifolds, and 
let $W$ be a submanifold of $P$. Let $\Gamma:N\times Z\to P$ be a mapping.
If $\Gamma$ is transverse to $W$, then there exists a subset $\Sigma$ of $Z$ with 
Lebesgue measure zero such that for any $p \in Z-\Sigma $, 
$\Gamma_p:N\to P$ is transverse to $W$, 
where $\Gamma_p(q)=\Gamma(q,p)$. 
\end{lemma}

\section{Proof of Theorem \ref{immersion}}\label{section 3}
Let $\{(U_\lambda ,\varphi _\lambda )\}_{\lambda \in \Lambda}$ be a coordinate neighborhood system of $N$. 
Let $\pi :J^1(N,\mathbb{R}^m)$$\to N\times \mathbb{R}^m$ be the natural projection defined by $\pi(j^1g(q))=(q,g(q))$. 
Let $\Phi _\lambda :\pi^{-1}(U_\lambda \times \mathbb{R}^m)\to \varphi _\lambda (U_\lambda)\times \mathbb{R}^m\times J^1(n,m)$ be 
the homeomorphism defined by 
\begin{eqnarray*} 
\Phi _\lambda \left(j^1g(q)\right)=\left(\varphi _\lambda (q),g(q),j^1(g\circ \varphi _\lambda ^{-1}\circ \widetilde{\varphi} _\lambda)(0)\right), 
\end{eqnarray*} 
where $\widetilde{\varphi} _\lambda : \mathbb{R}^n\to \mathbb{R}^n$ is the 
translation defined by $\widetilde{\varphi} _\lambda(0)=\varphi _\lambda (q)$. 
Then, $\{(\pi^{-1}(U_\lambda \times \mathbb{R}^m), 
\Phi _\lambda )\}_{\lambda \in \Lambda}$ 
is a coordinate neighborhood system of $J^1(N,\mathbb{R}^m)$. 
For any $k$ $(k=1,\ldots ,n)$, set 
\begin{eqnarray*} 
\Sigma ^k=\left\{j^1g(0)\in J^1(n,m)\mid {\rm dim\ Ker}Jg(0)=k\right\}.
\end{eqnarray*}
For any $k$ $(k=1,\ldots ,n)$, set 
\begin{eqnarray*}
\Sigma^k(N,\mathbb{R}^m)=\bigcup_{\lambda \in \Lambda}\Phi ^{-1}_\lambda \left(\varphi _\lambda (U_\lambda )\times \mathbb{R}^m \times \Sigma ^k\right). 
\end{eqnarray*}
Then, the set $\Sigma ^k(N,\mathbb{R}^m)$ is a subfiber-bundle of $J^1(N,\mathbb{R}^m)$ such that 
\begin{eqnarray*}
{\rm codim}\ \Sigma ^k(N,\mathbb{R}^m)&=&{\rm dim}\ J^1(N,\mathbb{R}^m)- 
{\rm dim}\ \Sigma ^k(N,\mathbb{R}^m) \\
&=&k(m-n+k).
\end{eqnarray*}
(for details on $\Sigma^k(N,\mathbb{R}^m)$, see for example 
\cite{GG}, pp.60--61). 

Now, let $\Gamma:N\times (\mathbb{R}^m)^m\to J^1(N,\mathbb{R}^m)$ 
be the mapping defined by 
\begin{eqnarray*}
\Gamma(q,p)=j^1(G_{(p,A)}\circ f)(q).
\end{eqnarray*}
We will show first that the mapping $\Gamma$ is transverse to the submanifold $\Sigma ^k(N,\mathbb{R}^m)$ 
for any $k$ $(k=1,\ldots,n)$. 
It is sufficient to show that if $\Gamma(\widetilde{q},\widetilde{p})\in \Sigma ^k(N,\mathbb{R}^m)$, then the following $(\ast)$ holds. 
\begin{align}
d\Gamma_{(\widetilde{q},\widetilde{p})}(T_{(\widetilde{q},\widetilde{p})}
(N\times (\mathbb{R}^m)^m))
+
T_{\Gamma(\widetilde{q},\widetilde{p})}\Sigma ^k(N,\mathbb{R}^m)
=
T_{\Gamma(\widetilde{q},\widetilde{p})}J^1(N,\mathbb{R}^m).\tag{$\ast $}
\end{align}
There exists a coordinate neighborhood 
$\left(U_{\widetilde{\lambda}}\times (\mathbb{R}^m)^m, \varphi_{\widetilde{\lambda}}\times id \right)$ 
containing the point $(\widetilde{q},\widetilde{p})$ of $N\times (\mathbb{R}^m)^m$, 
where $id$ is the identity mapping of $(\mathbb{R}^m)^m$ into $(\mathbb{R}^m)^m$, 
and the mapping $\varphi_{\widetilde{\lambda}}\times id : 
U_{\widetilde{\lambda}}\times (\mathbb{R}^m)^m \to \mathbb{R}^n\times (\mathbb{R}^m)^m$ is defined by 
$\left(\varphi_{\widetilde{\lambda}}\times id\right)(q,p)=
\left(\varphi_{\widetilde{\lambda}}(q), id(p)\right)$.
There exists a coordinate neighborhood 
$\left(\pi^{-1}(U_{\widetilde{\lambda}} \times \mathbb{R}^m), 
\Phi _{\widetilde{\lambda}} \right)$ 
containing the point $\Gamma (\widetilde{q},\widetilde{p})$ of 
$J^1(N,\mathbb{R}^m)$. 
Let $t=(t_1,\ldots ,t_n)\in \mathbb{R}^n$ be a local coordinate containing 
$\varphi_{\widetilde{\lambda}}(\widetilde{q})$. 
Then, the mapping $\Gamma$ is locally given by the following:
\begin{eqnarray*}
&{}&(\Phi _{\widetilde{\lambda}} \circ \Gamma \circ 
(\varphi_{\widetilde{\lambda}}\times id )^{-1})(t,p)\\
&=&(\Phi _{\widetilde{\lambda}} \circ \Gamma \circ 
(\varphi_{\widetilde{\lambda}}^{-1}\times id^{-1} ))(t,p)\\
&=&(\Phi _{\widetilde{\lambda}} \circ \Gamma )
(\varphi_{\widetilde{\lambda}}^{-1}(t), p )\\
&=&\Phi _{\widetilde{\lambda}}(\Gamma 
(\varphi_{\widetilde{\lambda}}^{-1}(t), p ))\\
&=&\Phi _{\widetilde{\lambda}}(j^1(G_{(p,A)}\circ f)(\varphi_{\widetilde{\lambda}}^{-1}(t)))\\
&=&(\Phi _{\widetilde{\lambda}}\circ j^1(G_{(p,A)}\circ f)\circ \varphi_{\widetilde{\lambda}}^{-1})(t)\\
&=&\left(t, (G_{(p,A)}\circ f \circ \varphi_{\widetilde{\lambda}}^{-1})(t), \right.  \\
&& \frac{\partial (G_1\circ f\circ \varphi_{\widetilde{\lambda}}^{-1})}
{\partial t_1}(t), 
\ldots ,
\frac{\partial (G_1\circ f\circ \varphi_{\widetilde{\lambda}}^{-1})}
{\partial t_n}(t),\\
&&\hspace{100pt}\cdots \cdots \cdots , \\
&& \left. \frac{\partial (G_m\circ f\circ \varphi_{\widetilde{\lambda}}^{-1})}
{\partial t_1}(t)
,\ldots ,
\frac{\partial (G_m\circ f\circ \varphi_{\widetilde{\lambda}}^{-1})}
{\partial t_n}(t)
\right)\\
&=&\left(t, (G_{(p,A)}\circ f \circ \varphi_{\widetilde{\lambda}}^{-1})(t), \right.  \\
&& 2\sum_{j=1}^m a_{1j}(\widetilde{f}_j(t)-p_{1j})
\frac{\partial \widetilde{f}_j}{\partial t_1}(t),\ldots ,2\sum_{j=1}^m a_{1j}
(\widetilde{f}_j(t)-p_{1j})
\frac{\partial \widetilde{f}_j}{\partial t_n}(t),\\
&&\hspace{100pt}\cdots \cdots \cdots , \\
&& \left. 2\sum_{j=1}^m a_{mj}(\widetilde{f}_j(t)-p_{mj})
\frac{\partial \widetilde{f}_j}{\partial t_1}(t),\ldots ,2\sum_{j=1}^m a_{mj}
(\widetilde{f}_j(t)-p_{mj})
\frac{\partial \widetilde{f}_j}{\partial t_n}(t)\right), 
\end{eqnarray*} 
where $p=(p_{11},\ldots, p_{1m},\ldots ,p_{m1},\ldots ,p_{mm})$, 
$f=(f_1,\ldots ,f_m)$, $G_{(p,A)}=(G_1,\ldots , G_m)$, and 
$\widetilde{f}_j= f_j\circ \varphi_{\widetilde{\lambda}}^{-1}$ $(1\leq j\leq m)$. 
The Jacobian matrix of the mapping $\Gamma$ at 
$(\widetilde{q},\widetilde{p})$ is the following:\\
\begin{eqnarray*}
J\Gamma_{(\widetilde{q},\widetilde{p})}=
\left(
\begin{array}{c|ccccccccccccccc}
E_n       &    0    &         \cdots          &   \cdots        &0    \\
\hline                &   \ast      &    \cdots    &  \cdots      & \ast  \\ 
      &  B_1 & & \bigzerol &\\
 \ast    &              &B_2     &       &  \\
    &    & \bigzerol &\ddots &      \\
      &    &  &  &  B_m\\
\end{array}
\right)_{(\varphi_{\widetilde{\lambda}}(\widetilde{q}),\widetilde{p})}, 
\end{eqnarray*}
where $E_n$ is the $n\times n$ unit matrix and $B_i$ $(1\leq i \leq m)$ is the 
following $n\times m$ matrix. 
\begin{eqnarray*}
B_i=\left(
\begin{array}{cccccc}
 -2a_{i1}\frac{\partial \widetilde{f}_1}{\partial t_1}(t) & \cdots & -2a_{im}\frac{\partial \widetilde{f}_m}{\partial t_1}(t)\\
    \vdots   & \ddots & \vdots        \\
 -2a_{i1}\frac{\partial \widetilde{f}_1}{\partial t_n}(t) & \cdots & -2a_{im}\frac{\partial \widetilde{f}_m}{\partial t_n}(t) \\
 \end{array}
\right)_{t=\varphi_{\widetilde{\lambda}}(\widetilde{q})}. 
\end{eqnarray*}
Since $\Sigma ^k(N,\mathbb{R}^m)$ is a subfiber-bundle of $J^1(N,\mathbb{R}^m)$ 
with fiber $\Sigma ^k$, 
in order to show $(\ast)$, 
it is clearly seen that the rank of the following matrix $C$ is $n+m+nm$. 
\begin{eqnarray*}
C=
\left(
\begin{array}{c|ccccccccccccccc}
E_{n+m}       &    \ast    &         \cdots          &   \cdots        &\ast   \\
\hline         &  B_1 & & \bigzerol &\\
\bigzerol  &              &B_2     &       &  \\
    &    & \bigzerol &\ddots &      \\
      &    &  &  &  B_m\\
\end{array}
\right)_{(\varphi_{\widetilde{\lambda}}(\widetilde{q}),\widetilde{p})}, 
\end{eqnarray*}
where $E_{n+m}$ is the $(n+m)\times (n+m)$ unit matrix. 
Notice that for any $i$ $(1\leq i \leq m^2)$, the $(n+m+i)$-th column vector of $C$ is 
the $(n+i)$-th column vector of $J\Gamma_{(\widetilde{q},\widetilde{p})}$. 
Let $Jf_{\widetilde{q}}$ be the Jacobian matrix of the mapping $f$ at $\widetilde{q}$. 
Since $a_{ij}\not=0$ for any $i$, $j$ $(1\leq i, j \leq m)$, 
there exists an $m\times m$ regular matrix $R$ such that 
$B_iR={}^t(Jf_{\widetilde{q}})$ 
for any $i$ $(1\leq i \leq m)$, where ${}^tX$ means the transposed matrix of $X$. 
Hence, there exists a 
$(n+m+m^2)\times (n+m+m^2)$ regular matrix $\widetilde{R}$ 
such that 
\begin{eqnarray*}
C\widetilde{R}=
\left(
\begin{array}{c|ccccccccccccccc}
E_{n+m}       &    \ast    &         \cdots          &   \cdots        &\ast   \\
\hline         &  {}^t(Jf_{\widetilde{q}})& & \bigzerol &\\
\bigzerol  &              &{}^t(Jf_{\widetilde{q}})   &       &  \\
    &    & \bigzerol &\ddots &      \\
      &    &  &  & {}^t(Jf_{\widetilde{q}})\\
\end{array}
\right)_{(\varphi_{\widetilde{\lambda}}(\widetilde{q}),\widetilde{p})}. 
\end{eqnarray*}
Since 
the mapping $f$ is an immersion $(n\leq m)$, we have that the rank of the matrix 
$C\widetilde{R}$ 
is $n+m+nm$. 
Therefore, the rank of the matrix $C$ must be $n+m+nm$. 
Hence, we have $(\ast)$. 
Thus, the mapping $\Gamma$ is transverse to the submanifold $\Sigma ^k(N,\mathbb{R}^m)$. 

By Lemma \ref{abra}, for any $k$ $(k=1,\ldots,n)$, 
there exists a subset $\widetilde{\Sigma}^k$ of $(\mathbb{R}^m)^m$ with 
Lebesgue measure zero such that 
for any $p\in (\mathbb{R}^m)^m-\widetilde{\Sigma}^k $, 
the mapping $\Gamma_p:N\to J^1(N,\mathbb{R}^m)$ is transverse to the submanifold 
$\Sigma ^k(N,\mathbb{R}^m)$. 
Set $\Sigma =\bigcup_{k=1}^n\widetilde{\Sigma }^k$. 
Notice that $\Sigma$ is a subset of $(\mathbb{R}^m)^m$ with Lebesgue measure zero. 
Then, for any $p\in (\mathbb{R}^m)^m-\Sigma $, 
the mapping $\Gamma_p:N\to J^1(N,\mathbb{R}^m)$ 
is transverse to the submanifold $\Sigma ^k(N,\mathbb{R}^m)$ for any $k$ $(k=1,\ldots,n)$. 

In order to show that for any $p\in (\mathbb{R}^m)^m-\Sigma $, 
the mapping $G_{(p,A)}\circ f:N\to \mathbb{R}^m$ is an immersion, 
it is sufficient to show that for any $p\in (\mathbb{R}^m)^m-\Sigma $, 
it follows that 
$\Gamma_p(N)\bigcap \bigcup_{k=1}^{n}\Sigma ^k(N,\mathbb{R}^m)=\emptyset $.

Suppose that there exists an element $p_0\in (\mathbb{R}^m)^m-\Sigma $ such that 
there exists an element $q_0\in N$ such that 
$\Gamma_{p_0}(q_0)\in \bigcup_{k=1}^{n}\Sigma ^k(N,\mathbb{R}^m)$. 
Then, there exists a natural number $k'$ $(1\leq k'\leq n)$ such that 
$\Gamma_{p_0}(q_0)\in \Sigma ^{k'}(N,\mathbb{R}^m)$. 
Since $\Gamma_{p_0}$ is 
transverse to $\Sigma ^{k'}(N,\mathbb{R}^m)$, we have the following:
\begin{eqnarray*}
d(\Gamma_{p_0})_{q_0}(T_{q_0}N)
+T_{\Gamma_{p_0}(q_0)}\Sigma ^{k'}(N,\mathbb{R}^m)
=T_{\Gamma_{p_0}(q_0)}J^1(N,\mathbb{R}^m).
\end{eqnarray*}
Hence, we have
\begin{eqnarray*}
{\rm dim}\ d(\Gamma_{p_0})_{q_0}(T_{q_0}N)
&\geq &{\rm dim}\ T_{\Gamma_{p_0}(q_0)}J^1(N,\mathbb{R}^m)-
{\rm dim}\ T_{\Gamma_{p_0}(q_0)}\Sigma ^{k'}(N,\mathbb{R}^m)\\
&=&{\rm codim}\ T_{\Gamma_{p_0}(q_0)}\Sigma ^{k'}(N,\mathbb{R}^m).
\end{eqnarray*}
Thus, we have $n\geq k'(m-n+k')$. 
This contradicts the assumptions $m\geq 2n$ and $k'\geq 1$. 
\hfill $\Box$
\section{Proof of Theorem \ref{injective}} \label{section 4}
Let $\Delta$ be the subset of $\mathbb{R}^{2m}$ defined by 
$\Delta =\{(y,y) \mid y\in \mathbb{R}^m\}$. 
It is clearly seen that $\Delta $ is a submanifold of $\mathbb{R}^{2m}$ 
such that 
\begin{eqnarray*}
{\rm codim}\ \Delta ={\rm dim}\ \mathbb{R}^{2m}- 
{\rm dim}\ \Delta  =m.
\end{eqnarray*}
Set $N^{(2)}=\{(q,q')\in N^2\mid q\not=q'\}$. 
Notice that $N^{(2)}$ is an open submanifold of $N^2$. 

Now, let $\Gamma : N^{(2)}\times (\mathbb{R}^m)^m\to \mathbb{R}^{2m}$ be the mapping defined by 
\begin{eqnarray*}
\Gamma(q,q',p)=\left( (G_{(p,A)}\circ f)(q), (G_{(p,A)}\circ f)(q')\right).
\end{eqnarray*}
We will show first that the mapping $\Gamma$ is transverse to the submanifold 
$\Delta $. 
It is sufficient to show that 
if $\Gamma(\widetilde{q},\widetilde{q}\,{'},\widetilde{p})\in \Delta $, 
then the following $(\ast \ast)$ holds. 
\begin{align}
d\Gamma_{(\widetilde{q},\widetilde{q}\,{'},\widetilde{p})}(T_{(\widetilde{q},\widetilde{q}\,{'},\widetilde{p})}(N^{(2)}\times (\mathbb{R}^m)^m))+
T_{\Gamma(\widetilde{q},\widetilde{q}\,{'},\widetilde{p})}\Delta
=T_{\Gamma(\widetilde{q},\widetilde{q}\,{'},\widetilde{p})}\mathbb{R}^{2m}.
\tag{$\ast \ast $}
\end{align}
Let $\{(U_\lambda ,\varphi _\lambda )\}_{\lambda \in \Lambda}$ be a coordinate  neighborhood system of $N$. 
There exists a coordinate neighborhood 
$\left(U_{\widetilde{\lambda}}\times U_{\widetilde{\lambda'}} \times (\mathbb{R}^m)^m, \varphi_{\widetilde{\lambda}}\times \varphi_{\widetilde{\lambda'}}\times id \right)$ 
containing the point $(\widetilde{q},\widetilde{q}\,{'},\widetilde{p})$ 
of $N^{(2)}\times (\mathbb{R}^m)^m$, 
where $id$ is the identity mapping of $(\mathbb{R}^m)^m$ into $(\mathbb{R}^m)^m$, 
and the mapping 
$\varphi_{\widetilde{\lambda}}\times \varphi_{\widetilde{\lambda'}}\times id : 
U_{\widetilde{\lambda}}\times U_{\widetilde{\lambda'}}\times (\mathbb{R}^m)^m \to \mathbb{R}^n\times \mathbb{R}^n\times (\mathbb{R}^m)^m$ is defined by 
$\left(\varphi_{\widetilde{\lambda}}\times \varphi_{\widetilde{\lambda'}}\times id\right)(q,q',p)=
\left(\varphi_{\widetilde{\lambda}}(q), \varphi_{\widetilde{\lambda '}}(q'), id(p)\right)$.
Let $t=(t_1,\ldots ,t_n)$ be a local coordinate 
containing $\varphi_{\widetilde{\lambda}}(\widetilde{q})$, 
and let $t'=(t_1'\ldots ,t_n')$ be a local coordinate 
containing $\varphi_{\widetilde{\lambda'}}(\widetilde{q}\,{'})$. 
Then, the mapping $\Gamma$ is locally given by the following:
\begin{eqnarray*}
&&\Gamma \circ \left(\varphi_{\widetilde{\lambda}}\times \varphi_{\widetilde{\lambda'}}\times id\right)^{-1}(t,t',p)
\\
&=&\Gamma \circ \left(\varphi_{\widetilde{\lambda}}^{-1}\times \varphi_{\widetilde{\lambda'}}^{-1}\times id^{-1}\right)(t,t',p)
\\
&=&\Gamma \left(\varphi_{\widetilde{\lambda}}^{-1}(t),  
\varphi_{\widetilde{\lambda'}}^{-1}(t'), p\right)
\\
&=&\left( (G_{(p,A)}\circ f\circ \varphi_{\widetilde{\lambda}}^{-1})(t), 
(G_{(p,A)}\circ f\circ \varphi_{\widetilde{\lambda'}}^{-1})(t')\right)
\\
&=&\left(\sum_{j=1}^m a_{1j}(\widetilde{f}_j(t)-p_{1j})^2,\ldots ,
\sum_{j=1}^m a_{mj}(\widetilde{f}_j(t)-p_{mj})^2,\right.\\
&& \left. \sum_{j=1}^m a_{1j}(\widetilde{f}'_j(t')-p_{1j})^2,\ldots ,
\sum_{j=1}^m a_{mj}(\widetilde{f}'_j(t')-p_{mj})^2\right), 
\end{eqnarray*}  
where $p=(p_{11},\ldots, p_{1m},\ldots ,p_{m1},\ldots ,p_{mm})$, 
$f=(f_1,\ldots ,f_m)$, 
$\widetilde{f}_j= f_j\circ \varphi_{\widetilde{\lambda}}^{-1}$, 
and $\widetilde{f}'_j= f_j\circ \varphi_{\widetilde{\lambda}'}^{-1}$ 
$(1\leq j\leq m)$. 
The Jacobian matrix of the mapping $\Gamma$ at 
$(\widetilde{q},\widetilde{q}\,{'},\widetilde{p})$ is the following: 
\begin{eqnarray*}
J\Gamma_{(\widetilde{q},\widetilde{q}\,{'},\widetilde{p})}=
\left(
\begin{array}{c|cccccccccccccc}
&{\bf b}_{1}& & &\bigzerol \\
&&{\bf b}_{2}&& \\
&\bigzerol &   &  \ddots   &    \\ 
\ast&&   &         &{\bf b}_{m}\\ 
 &{\bf b'}_{1}&  & &   \bigzerol       \\
&& {\bf b'}_{2} & &\\
&\bigzerol&   &  \ddots   &    \\ 
&  & &&{\bf b'}_{m}     \\
\end{array}
\right)_{(\varphi_{\widetilde{\lambda}}(\widetilde{q}), \varphi_{\widetilde{\lambda}'}(\widetilde{q}\,{'}), \widetilde{p})},
\end{eqnarray*}
where 
\begin{eqnarray*}
{\bf b}_{i}=-2\left(
a_{i1}(\widetilde{f}_1(t)-p_{i1}),\ldots ,a_{im}(\widetilde{f}_m(t)-p_{im})
\right), \\
{\bf b'}_{i}=-2\left(
a_{i1}(\widetilde{f}'_1(t')-p_{i1}),\ldots ,a_{im}(\widetilde{f}'_m(t')-p_{im})
\right). 
\end{eqnarray*}
By seeing the construction of 
$T_{\Gamma(\widetilde{q},\widetilde{q}\,{'},\widetilde{p})}\Delta$, 
in order to show $(\ast \ast)$, it is sufficient to 
show that the rank of the following matrix $D$ is $2m$. 
\begin{eqnarray*}
D=
\left(
\begin{array}{c|cccccccccccccc}
&{\bf b}_{1}& & &\bigzerol \\
E_m&&{\bf b}_{2}&& \\
&\bigzerol &   &  \ddots   &    \\ 
&&   &         &{\bf b}_{m}\\ 
 \hline  &{\bf b'}_{1}&  & &   \bigzerol       \\
E_m&& {\bf b'}_{2} & &\\
&\bigzerol&   &  \ddots   &    \\ 
&  & &&{\bf b'}_{m}     \\
\end{array}
\right)_{(\varphi_{\widetilde{\lambda}}(\widetilde{q}), \varphi_{\widetilde{\lambda}'}(\widetilde{q}\,{'}), \widetilde{p})},
\end{eqnarray*}
where $E_m$ is the $m\times m$ unit matrix. 
Notice that for any $i$ $(1\leq i \leq m^2)$, 
the $(m+i)$-th column vector of $D$ is 
the $(2n+i)$-th column vector of 
$J\Gamma_{(\widetilde{q},\widetilde{q}\,{'},\widetilde{p})}$. 
By $a_{ij}\not=0$, there exist an $(m+m^2)\times (m+m^2)$ regular matrix 
$Q_1$ such that the following holds:
\\
\scalebox{0.85}{\begin{minipage}{\linewidth}
\begin{eqnarray*}
&&DQ_1=
\\&&
\left(
\begin{array}{c|cccccccccccccc}
&\widetilde{f}_1(t)-p_{11}& \cdots &\widetilde{f}_m(t)-p_{1m} \\
E_m && 0 &            & \ddots &&0\\
&& & &&\widetilde{f}_1(t)-p_{m1}&\cdots&\widetilde{f}_m(t)-p_{mm} \\
\hline \\ \\[-7mm]
&\widetilde{f}'_1(t')-p_{11}& \cdots &\widetilde{f}'_m(t')-p_{1m} \\
E_m && 0 &            & \ddots &&0\\
&& & &&\widetilde{f}'_1(t')-p_{m1}&\cdots&\widetilde{f}'_m(t')-p_{mm} \\
\end{array}
\right)_{(t,t',p)},
\end{eqnarray*}\end{minipage}}\\
where $(t,t',p)=(\varphi_{\widetilde{\lambda}}(\widetilde{q}), \varphi_{\widetilde{\lambda}'}(\widetilde{q}\,{'}), \widetilde{p})$. 
It is clearly seen that there exist a $2m\times 2m$ regular matrix $Q_2$ 
and an $(m+m^2)\times (m+m^2)$ regular matrix $Q_3$ such that the following 
holds: \\
\scalebox{0.8}{\begin{minipage}{\linewidth}
\begin{eqnarray*}
&&Q_2DQ_1Q_3=
\\&&
\left(
\begin{array}{c|cccccccccccccc}
&& &   & & \\
E_m&& &&\bigzerol &&    \\
&&  &      &  &&  &  \\ 
\hline \\ \\[-7mm]
&\widetilde{f}'_1(t')-\widetilde{f}_1(t)& \cdots &\widetilde{f}'_m(t')-\widetilde{f}_m(t)  \\
0 && 0 &            & \ddots &&0\\
&& & &&\widetilde{f}'_1(t')-\widetilde{f}_1(t)&\cdots&\widetilde{f}'_m(t')-\widetilde{f}_m(t) \\
\end{array}
\right)_{(t,t',p)},
\end{eqnarray*}\end{minipage}}\\
where $(t,t',p)=(\varphi_{\widetilde{\lambda}}(\widetilde{q}), \varphi_{\widetilde{\lambda}'}(\widetilde{q}\,{'}), \widetilde{p})$. 
Since $f$ is injective, there exists a natural number $j$ $(1\leq j\leq m)$ such that 
$\widetilde{f}'_{j}(t')-\widetilde{f}_{j}(t)\not=0$. 
Hence, we have that the rank of $Q_2DQ_1Q_3$ is $2m$. 
Therefore, the rank of the matrix $D$ must be $2m$. 
Hence, we have $(\ast \ast)$. 
Thus, the mapping $\Gamma$ is transverse to the submanifold $\Delta$. 

By Lemma \ref{abra}, 
there exists a subset $\Sigma$ of $(\mathbb{R}^m)^m$ with 
Lebesgue measure zero such that 
for any $p\in (\mathbb{R}^m)^m-\Sigma $, 
the mapping $\Gamma_p:N^{(2)}\to \mathbb{R}^{2m}$ is transverse to 
the submanifold $\Delta$. 

In order to prove that for any $p\in (\mathbb{R}^m)^m-\Sigma$, 
the mapping $G_{(p,A)}\circ f$ is injective, 
it is sufficient to show that for any $p\in (\mathbb{R}^m)^m-\Sigma$, 
it follows that $\Gamma_{p}(N^{(2)})\cap \Delta =\emptyset$. 
Suppose that there exists an element $p_0\in (\mathbb{R}^m)^m-\Sigma$ such that 
there exists an element $(q_0,q_0')\in N^{(2)}$ such that 
$\Gamma_{p_0}(q_0,q_0')\in \Delta $. 
Since $\Gamma_{p_0}$ is 
transverse to $\Delta $, we have the following:
\begin{eqnarray*}
d(\Gamma_{p_0})_{(q_0,q_0')}(T_{(q_0,q_0')}N^{(2)})
+T_{\Gamma_{p_0}(q_0,q_0')}\Delta
=T_{\Gamma_{p_0}(q_0,q_0')}\mathbb{R}^{2m}.
\end{eqnarray*}
Hence, we have 
\begin{eqnarray*}
{\rm dim}\ d(\Gamma_{p_0})_{(q_0,q_0')}(T_{(q_0,q_0')}N^{(2)})
&\geq &{\rm dim}\ T_{\Gamma_{p_0}(q_0,q_0')}\mathbb{R}^{2m}-
{\rm dim}\ T_{\Gamma_{p_0}(q_0,q_0')}\Delta 
\\
&=&{\rm codim}\ T_{\Gamma_{p_0}(q_0,q'_0)}\Delta.
\end{eqnarray*}
Thus, we have $2n\geq m$. 
This contradicts the assumption $m\geq 2n+1$. 
\hfill $\Box$

\section{Appendix} \label{section 5}
The proofs of the following Lemma \ref{singularities} and Lemma \ref{not injective}
are given in Subsection \ref{singularities} and Subsection \ref{not injective}, 
respectively. 
\begin{lemma}\label{singularities}
Any generalized distance-squared mapping of equidimensional 
case $G_{(p,A)}:\mathbb{R}^m\to \mathbb{R}^m$ has singularities. 
\end{lemma}
\begin{lemma}\label{not injective}
Any generalized distance-squared mapping of equidimensional 
case $G_{(p,A)}:\mathbb{R}^m\to \mathbb{R}^m$ is not injective. 
\end{lemma}
\subsection{Proof of Lemma \ref{singularities}}\label{singularities}
Let $J(G_{(p,A)})_x$ be the Jacobian matrix of the mapping $G_{(p,A)}$ at $x$. 
\begin{eqnarray*}
J(G_{(p,A)})_x=
2\left(
\begin{array}{ccccccccccccccc}
a_{11}(x_1-p_{11})&\cdots  &a_{1m}(x_m-p_{1m}) \\
\vdots  & \ddots     &\vdots  \\ 
a_{m1}(x_1-p_{m1})&\cdots &a_{mm}(x_m-p_{mm})   \\
\end{array}
\right)_x. 
\end{eqnarray*}
If $x=p_i$ $(1\leq i\leq m)$, then we have that ${\rm rank} J(G_{(p,A)})_{p_i}\leq m-1$. 
\hfill $\Box$
\subsection{Proof of Lemma \ref{not injective}}\label{not injective}
Set $G_{(p,A)}=(G_1, \ldots, G_m)$.    It is clear that 
$G_{(p, A)}^{-1}(\{0\}\times \mathbb{R}^{m-1})=G_1^{-1}(0)$.    
Since $G_1$ has the form 
$G_1(x)=\sum_{j=1}^m a_{1j}(x_j-p_{1j})^2$ $(a_{1j}\ne 0)$, 
it is easy to see that $G_{1}^{-1}(0)=\{p_1\}$ or 
$G_{1}^{-1}(0)-\{p_1\}$ is homotopy equivalent to 
$S^k\times S^{m-2-k}$ where $k$ is an integer such that 
$0\le k\le m-2$.    
Hence, it follows that the set-germ 
$(\{0\}\times \mathbb{R}^{m-1}, G_{(p,A)}(p_1))$ is not homeomorphic 
to the set-germ $(G_1^{-1}(0), p_1)$.  
\par 
On the other hand, suppose that $G_{(p,A)}$ is injective.   
Then, by the invariance of domain theorem (\cite{spanier}),  
$G_{(p,A)}^{-1}: G_{(p,A)}(\mathbb{R}^m)\to \mathbb{R}^m$ 
must be a homeomorphism.    
It follows that the set-germ 
$(\{0\}\times \mathbb{R}^{m-1}, G_{(p,A)}(p_1))$ is homeomorphic 
to the set-germ $(G_1^{-1}(0), p_1)$, which is a contradiction.    
Therefore, $G_{(p,A)}$ is not injective.   

\hfill $\Box$
\section*{Acknowledgements}
The first author is supported by JSPS KAKENHI Grant Number 16J06911.

\end{document}